\documentclass{article}
\usepackage{graphicx}
\usepackage{amssymb}
\usepackage{amsmath}
\usepackage{amscd}
\usepackage{amsthm} 

\usepackage[utf8]{inputenc}
\usepackage{amsfonts}
\usepackage{amstext}
\usepackage{color}
\newtheorem{theorem}{Theorem}
\newtheorem{lemma}{Lemma}

 \newtheorem{remark}{Remark}\newtheorem{proposition}{Proposition}
 
\newcommand{\R}{{\mathbb R}}  \newcommand{\Z}{{\mathbb Z}} \newcommand{\N}{{\mathbb N}}

\newcommand{\Ee}{{\mathbb E}}\newcommand{\al}{\alpha}

\providecommand{\keywords}[1]
{
  \small	
  \textbf{\textit{Keywords: }} #1
}

\begin{document}
\author{Aleksei Kulikov \quad Alexander~Ulanovskii \quad Ilya~Zlotnikov }

\title{Completeness  of Certain Exponential Systems and Zeros of Lacunary Polynomials}
\maketitle
\begin{abstract}

Let $\Gamma$ be a subset of $\{0,1,2,...\}$. We show that if  $\Gamma$ has `gaps' then the completeness and frame properties of the system $\{t^ke^{2\pi i nt}: n\in\Z,k\in\Gamma\}$   differ from those of  the classical exponential systems. This phenomenon is closely connected with the existence of certain uniqueness  sets for lacunary polynomials.
\end{abstract}
{\keywords{completeness, frame, totally positive matrix, generalized Vandermonde matrix, uniqueness set, lacunary polynomials}}

 \date{}\maketitle
 
\section{Introduction} 
Let $\Lambda$ be a {\it separated} set of real numbers. Denote by $$E(\Lambda):=\{e^{2\pi i\lambda t},\lambda\in\Lambda\}$$ the corresponding exponential system.

Approximation and representation properties of exponential systems in different function spaces is a classical subject of investigation. In particular,  the completeness and the frame problems of $E(\Lambda)$ for the space $L^2(a,b)$ can be stated as follows: Determine if

\begin{enumerate}
    \item[(a)] ({\it Completeness property of $E(\Lambda)$}) every function $F$ in $L^2(a,b)$ can be  approximated arbitrarily well in $L^2$-norm by  finite linear combinations of exponential functions from $E(\Lambda)$;

    \item[(b)]({\it Frame property of $E(\Lambda)$}) there exist
two positive constants $A$ and $B$ such that for every $F\in L^2(a,b)$ we have 
$$A\|F\|^2_{2}\leq \sum_{\lambda\in\Lambda}|\langle F,e^{2\pi i\lambda t}\rangle|^2\leq B\|F\|^2_2,$$where $\langle\cdot,\cdot\rangle$ is the usual inner product in $L^2(a,b)$. 

\end{enumerate}

\medskip
Note that the notion of frame is very important and can be defined in similar manner for an  arbitrary system  of elements $E=\{e_\lambda\}$   in a Hilbert space $H$.
If $E$ is a frame in $H$, then every element $f$ from $H$
admits a (maybe, non-unique) representation 
$$f=\sum_{e_\lambda\in E}c_\lambda e_\lambda,$$
for some $l^2-$sequence of complex numbers $c_\lambda$ (see e.g. \cite{cr}).

It is easy to check that the completeness property of $E(\Lambda)$ is translation-invariant: If $E(\Lambda)$ is complete in $L^2(a,b)$, then it is complete in $L^2(a+c,b+c)$, for every $c\in\R$. 
As a `measure of completeness', one may introduce the so-called {\it completeness radius} of $E(\Lambda)$:
$$CR(\Lambda)=\sup\{a\geq 0: E(\Lambda) \mbox{ is complete in } L^2(-a,a)\}.$$Similarly, the frame property of $E(\Lambda)$ is also translation-invariant, and one may introduce the {\it frame radius} as 
$$FR(\Lambda)=\sup\{a\geq 0: E(\Lambda) \mbox{ is a frame in } L^2(-a,a)\}.$$

Both  radii above can be expressed  in terms of certain densities:

\medskip
(A) The celebrated Beurling--Malliavin theorem \cite{BM} states that $CR(\Lambda)=D^\ast(\Lambda)$. Here $D^\ast$ is the so-called upper (or external) Beurling--Malliavin density.

(B) It follows from the classical `Beurling Sampling Theorem' \cite{B} (see also a detailed discussion  in \cite{ou1}) that $FR(\Lambda)=D^-(\Lambda)$, where $\Lambda$ is a separated (also called uniformly discrete) set and  $D^-(\Lambda)$ is the lower uniform density of $\Lambda$.

We refer the reader to  \cite{OS} or \cite{S} for a complete description of exponential frames for the space $L^2(a,b)$. It is not given in terms of a density of $\Lambda$.

\medskip
Observe that the proofs of (A) and (B)   use techniques from the complex analysis.

The density $D^\ast$ can be defined and the Beurling--Malliavin formula for the completeness radius  remains valid 
for the  {\it multisets} $(\Lambda, \Gamma(\lambda))$, where $\Lambda\subset\R$ and $\Gamma(\lambda)=\{0,...,n(\lambda)-1\}$, i.e. for the systems 
\begin{equation}\label{es}E(\Lambda, \Gamma(\lambda)):=\{t^k e^{2\pi i\lambda t}: \lambda\in\Lambda, t=0,...,n(\lambda)-1 \}.\end{equation}
Here $n(\lambda)$ is the multiplicity (number of occurrences) of the element $\lambda\in\Lambda$. 
The same is true for the frame radius, see \cite{Gr}.
In particular, if  $\Lambda=\Z$  and $\Gamma(\lambda)=\Gamma_N:=\{0,...,N-1\}$,  $\lambda\in\Lambda$, then one has 
\begin{equation}\label{o}CR(\Z,\Gamma_N)=FR(\Z,\Gamma_N)=N/2=\#\Gamma_N/2,\end{equation}where $\#\Gamma$ is the number of elements of $\Gamma$,
$CR(\Z,\Gamma_N)$ and $FR(\Z,\Gamma_N)$ are the completeness and frame radius of $E(\Z,\Gamma_N)$, respectively.

One may consider the completeness property of systems in (\ref{es})  in  $L^p(a,b)$ and $C([a,b])$. For each of these spaces, the completeness property is translation-invariant. Clearly, the completeness in $C([-a,a])$ implies the completeness in $L^p(-a,a)$ for every $1\leq p <\infty.$
Observe that if  $E(\Lambda,\Gamma(\Lambda))$ is not complete in $C([-a,a])$,  its  {\it deficiency}  in $C([-a,a])$ is at most $1$, i.e. by adding to the system an  exponential function $e^{2\pi i at },a\not\in\Lambda,$ the new lager system becomes complete in $C([-a,a])$ (see e.g. discussion in \cite{r}).
It easily follows that   every system in (\ref{es}) has the same completeness radius   for every space considered above.

\section{Statement of Problem and Results}
Let us now introduce somewhat  more general  systems. Assume that $\Lambda\subset\R$ is a discrete set and that to every $\lambda\in\Lambda$ there corresponds  a finite or infinite set $\Gamma(\lambda)\subset \N_0:=\{0,1,2,3,...\}$. Set
$$E(\Lambda, \Gamma(\lambda))=\{t^\gamma e^{2\pi i \lambda t}: \lambda\in\Lambda, \gamma\in \Gamma(\lambda)\}.$$ 

Inspired by a recent work of H. Hedenmalm \cite{Hedenmalm}, we ask:
What are the completeness and frame properties of $E(\Lambda,\Gamma(\lambda))$? 
In this note we restrict ourselves to the case  $\Lambda=\Z$ and $\Gamma(n)=\Gamma\subset\N_0, n\in\Z,$ is a fixed set. That is, we will consider the completeness and frame properties of the system 
$$E(\Z,\Gamma):=\{t^\gamma e^{2\pi i n t}: n\in\Z, \gamma\in \Gamma\}, \quad \Gamma\subset\N_0.$$

Let us now introduce the formal analogues of the completeness and frame radius:
$$
CR(\Z,\Gamma):=\sup\{a\geq0: E(\Z,\Gamma)\mbox{ is complete in } L^2(-a,a)\},$$$$ FR(\Z,\Gamma):=\sup\{a\geq0: E(\Z,\Gamma)\mbox{ is a frame in } L^2(-a,a)\}.
$$
We also define the completeness radius   $CR_C(\Z,\Gamma)$          in the spaces of continuous functions:
$$
CR_C(\Z,\Gamma):=\sup\{a\geq0: E(\Z,\Gamma)\mbox{ is complete in } C([-a,a])\}.$$

In what follows, to exclude trivial remarks, we will always assume that $0\in\Gamma$. 

Set
    $$\Gamma_{even} = \Gamma \cap 2 \Z \quad \text{and} \quad \Gamma_{odd} = \Gamma \cap (2\Z + 1),$$and introduce the following number
    $$
r(\Gamma) := \begin{cases}
\# \Gamma_{odd} + \frac{1}{2} , \text{   if   } \# \Gamma_{odd} < \# \Gamma_{even}, \\[.2cm]
   \# \Gamma_{even} , \text{   if   } \# \Gamma_{odd} \ge \# \Gamma_{even}.\\[.2cm]  
\end{cases}
$$

Observe that $r(\Gamma)<\#\Gamma/2$ unless  $\#\Gamma_{even} = \#\Gamma_{odd}$ or $\#\Gamma_{even} = \#\Gamma_{odd}+1$.

It turns out  that  the completeness and frame properties of $E(\Z,\Gamma)$ may differ from the ones for the systems considered above. In particular, we have
\begin{theorem}\label{t1} Given any finite or infinite set $\Gamma\subset\N_0$ satisfying $0\in\Gamma.$ Then

{\rm (i)} $CR(\Z,\Gamma)=\#\Gamma/2$;

{\rm(ii)} $CR_C(\Z,\Gamma)=FR(\Z,\Gamma)=r(\Gamma).$

    \end{theorem}
  
 Below we  prove  more precise results. 

Theorem \ref{t1} shows that property   (\ref{o}) is no longer true for the systems  $E(\Z,\Gamma)$.
 
 The proof of part (i) uses mainly basic linear algebra. 
 We will see that the completeness property of $E(\Z,\Gamma)$ in $L^2(a,b)$ is translation invariant, and so $CR(\Z,\Gamma)$ still can be viewed as  a `measure of completeness' of $E(\Z,\Gamma)$.
 
 On the other hand, neither the frame property in $L^2(a,b)$ nor the completeness property in  $C([a,b])$ is translation invariant in the sense that both of them  depend  on the length of the interval $(a,b)$ and  also on its position. This phenomenon is intimately connected with the solvability of certain systems of linear equations and also with the existence of certain uniqueness sets for lacunary polynomials, see Theorem \ref{t2} below.

Given any finite set $M\subset\N_0$, let $P(M)$ denote the set of real polynomials with exponents in $M$: $$P(M):=\{P(x)=\sum_{m_j\in M}c_j x^{m_j}: c_j\in\R\}.$$

If $M\subset\N_0$ consists of  $n$ elements (shortly, $\#M=n$), then clearly no set  $X\subset\R$ satisfying $\# X\leq n-1$ is a uniqueness set for $P(M)$, i.e. there is a non-trivial polynomial $P\in P(M)$ which vanishes on $X$. This is no longer true if  $\# X=n$. Moreover,  there exist real uniqueness sets $X$,   $\# X=n$, that are uniqueness sets for every space $P(M),\# M=n$. Indeed, by Descartes' rule of signs, each $P\in P(M)$ may have at most $n-1$ distinct positive zeros, and so every set of $n$  positive points is a uniqueness set for $P(M)$.  Here we present  a less trivial example of such sets. Given $N$ distinct real numbers $t_1,\dots, t_N,$ set
\begin{equation}\label{ss}
S(t_1, \ldots , t_N) := \{ (-1)^kt_k\}_{k = 1}^N.
\end{equation}

\begin{theorem}\label{t2} Assume $0 < t_1<t_2<\dots <t_N$.  Then both sets $\pm S(t_1,\dots,t_N)$  are
uniqueness sets for every space $P(M),M\subset\N_0, \# M= N.$
\end{theorem}

The rest of the paper is organized as follows: In Section 3 several auxiliary results are proved. Theorem \ref{t2} is proved in Section 4.
We consider the completeness property of $E(\Z,\Gamma)$ in $L^2(a,b)$ and in $C([a,b])$ in Sections 5 and 6, respectively. Finally, in Section 7 we consider the frame property of $E(\Z,\Gamma) $ and also present some remarks. 

\section{Auxiliary Lemmas} 

Given $N \in \N, {\bf x} = \{x_0, \dots, x_{N-1}\} \subset \R$, and $ \Gamma = \{\gamma_0, \gamma_1, \dots, \gamma_{N-1}\} \subset \N$ we denote by $V({\bf x},\Gamma)$ a {\it generalized $N\times N$ Vandermonde matrix},
\begin{equation}
V({\bf x};\Gamma):= \begin{pmatrix}
x_0^{\gamma_0}& x_1^{\gamma_0}& x_2^{\gamma_0}& \dots&  x_{N-1}^{\gamma_0}&  \\[0.1cm]
x_0^{\gamma_1}& x_1^{\gamma_1}& x_2^{\gamma_1}& \dots&  x_{N-1}^{\gamma_1}&  \\[0.1cm]
\dots& \dots& \dots& \dots&  \dots& \\[0.1cm]
x_0^{\gamma_{N-1}}& x_1^{\gamma_{N-1}}& x_2^{\gamma_{N-1}}& \dots&  x_{N-1}^{\gamma_{N-1}}&  \\[0.1cm]
\end{pmatrix}.
\end{equation}
We will usually assume that $0\in\Gamma$. Note that if $\Gamma = \{0,1, \dots, N-1\}$, then the matrix $V({\bf x};\Gamma)$ is a standard Vandermonde matrix, and it is easy to compute its determinant and establish whenever it is invertible or not. However, if $\Gamma$ has gaps, the situation is more complicated.  In the case when $x_i > 0$ for all $i = 0, \dots, n-1$, one may use the following result from the theory of {\it totally positive matrices}, see e.g. \cite{Karlin} and \cite{Pinkus}.
\begin{proposition}{\rm(see \cite{Pinkus}, section~4.2)}\label{pos_matr_prop}
If $0 < x_0 < x_1 < \dots < x_N$ and $\gamma_0<\gamma_1 < \gamma_2 < \dots < \gamma_N$, then $V({\bf x}; \Gamma)$ is a totally positive matrix. In particular, it is  invertible.
\end{proposition}

This statement is no longer true  if ${\bf x}$ contains both positive and negative coordinates. 

We will be interested in a particular case where ${\bf x}=(s,s+1,...,s+N-1)$ for some $s\in\R.$ 
Consider the problem: Describe the set of points $s\in\R$
such that the matrix $V((s,\dots,s+N-1);\Gamma)$ is invertible {\it for every} $\Gamma\subset \N_0,\#\Gamma=N$.

\begin{lemma}\label{l0} $V((x_0, x_1,\ldots,x_{N-1});\Gamma)$ is not invertible if and only if there exists a polynomial $P\in P(\Gamma)$ which vanishes on  the set $\{x_0,x_1,...,x_{N-1}\}$.
\end{lemma}

\begin{proof}
Write $\Gamma=\{\gamma_0,\gamma_1,\dots,\gamma_{N-1}\}$. The matrix $V((x_0, x_1, \ldots , x_{N-1}); \Gamma)$ is not invertible if and only if its transpose is not. The latter means that there is a
 non-zero vector ${\bf a} = (a_0, \ldots , a_{N-1})$ satisfying  $V((x_0,x_1,\dots,x_{N-1}); \Gamma)^T{\bf a}^T = 0$. This means  that the polynomial $\sum_{j = 0}^{N-1} a_j x^{\gamma_j}$ vanishes at the points $x_0,\dots,x_{N-1}$.\end{proof}

\begin{lemma}\label{c1} Given $N\geq 2,$ the matrix $V((s,\dots,s+N-1);\Gamma)$ is invertible for every $\Gamma\subset\N_0,\#\Gamma=N,$ $0\in \Gamma$, if 

{\rm(i)} $s\geq 0;$

{\rm(ii)} $s\in (-N/2,-N/2+1)\setminus (1/2)\Z.$ 

\end{lemma}

Part (i) is a direct consequence of 
Proposition \ref{pos_matr_prop}.  

Part (ii) follows from Lemma~\ref{l0},  Theorem \ref{t2}, and the observation that for every $s\in (-N/2,-N/2+1)$ such that  $s$ does not equal $k/2$ for some $k\in\Z,$ the set  $\{s,\dots, s+N-1\}$ can be written as $\pm S$, where $S$ is defined  in~(\ref{ss}).

Clearly, by Lemma \ref{c1},  the determinant of $V((s,\dots,s+N-1);\Gamma)$ is a non-trivial polynomial of $s$. Hence, for every fixed $\Gamma$, this matrix is invertible for every $s$ outside a  finite number of points. 

In what follows, by measure we mean a finite, complex Borel measure on $\R$.

Given a  measure $\mu$, as usual we denote by $\hat\mu$ its Fourier-Stieltjes transform
$$\hat\mu(x)=\int\limits_\R e^{-2\pi ixt}\,d\mu(t).$$We also denote by $\delta_x$ the $\delta$-measure concentrated at the point  $x$.

\begin{lemma}\label{l1}
Let $\mu$ be a  measure supported by an interval $[\al,\al+1]$. The following are equivalent:

{\rm (i)}  $\hat\mu$ vanishes on $\Z;$

{\rm(ii)}  $\mu=A(\delta_\al-\delta_{\al+1})$, for some $A\in\mathbb{C}$.
\end{lemma}

\begin{proof}
We present a proof of  (i) $\Rightarrow $ (ii). The converse implication is trivial.

Since supp$\,\mu\subset[\al,\al+1]$, it is easy to see that the entire function 
$$f(z):=e^{2\pi i (\al+1/2)z}\hat\mu(z)$$
satisfies
\begin{equation}\label{f1}|f(x+iy)|\leq Ce^{\pi|y|},\quad x,y\in\R,\end{equation}
with some constant $C.$ Since $f$ vanishes on $\Z$, the function
$g(z):=f(z)/(\sin\pi z)$ is also entire. Clearly, there is a positive constant $B$ such that 
$$|\sin(\pi(x+iy))|\geq Be^{\pi|y|},\quad \mbox{for all } x,y\in\R, \ \inf_{n\in\Z}|x+iy-n|\geq 1/4.$$This, (\ref{f1}) and  the maximum modulus principle imply that   $g(z)$ is bounded in $\mathbb C$. hence, $g$ is a constant function, from which the lemma follows.   
\end{proof}

Let us now  consider measures $\mu$ that are "orthogonal" to $E(\Z,\Gamma)$:
\begin{equation}\label{orth}
\int\limits_\R t^\gamma e^{-2\pi i n t}\,d\mu(t)=0,\quad \mbox{for all } \gamma\in\Gamma, n\in\Z.
\end{equation}

\begin{lemma}\label{l2}
Assume that $\Gamma\subset\N_0,\#\Gamma=N,0\in \Gamma$, and that a measure $\mu$ is concentrated on  $[\al,\al+N]$. If $\mu$ satisfies {\rm (\ref{orth})}, then
there is a finite set $S\subset (\al,\al+1)$ and measures $\mu_s,s\in S,$ and $\nu$ such that

{\rm (i)} $\mu=\sum\limits_{s\in S}\mu_s+\nu$;

{\rm (ii)} $\nu$ and $\mu_s,s\in S,$ satisfy {\rm (\ref{orth})};

{\rm (iii)} The representations are true:
\begin{equation}\label{rep}
  d\nu=\sum_{j=1}^{N+1}a_j\delta_{\al+j-1}, \ \ d\mu_s= \sum_{j=1}^N c_{s,j}\delta_{s+j-1},\quad s\in S, \, c_{s,j} \in \R,\, a_j \in \R.
\end{equation}\end{lemma}

Note that $\mu_s$ satisfies {\rm (\ref{orth})} if and only if 
\begin{equation}\label{coef}\sum_{j=1}^N(s+j-1)^\gamma c_{s,j}=0,\quad \mbox{for every } \gamma\in\Gamma,\ s\in S.
\end{equation}
A similar observation is true for the measure $\nu$.

\begin{proof}[Proof of Lemma \ref{l2}]
 Clearly, $\mu$ admits a  unique representation
\begin{equation}\label{mu-decomp}
d\mu(x)=\sum_{j=1}^{N}d\mu_j(x-j+1),\end{equation}
where each $\mu_j$ is a measure supported by $[\al,\al+1)$ for $j=1,\dots, N-1,$ and supp$\,\mu_N=[\al,\al+1]$.
Then (\ref{orth}) is equivalent to
$$
\int_{[\al,\al+1]}e^{-2\pi i n t}\sum_{j=1}^N(t+j-1)^\gamma\,d\mu_j(t)=0,\quad \text{for every   } \gamma\in\Gamma, n\in\Z.
$$

It follows from Lemma \ref{l1} that $\mu_j$ satisfy the system of $N$ equations
\begin{equation}\label{seq}
\sum_{j=1}^N (t+j-1)^\gamma d\mu_j(t)=C_\gamma (\delta_\al-\delta_{\al+1}),\quad \text{for every   } \gamma\in\Gamma.
\end{equation}The corresponding matrix  on the left hand-side is $V((t,\dots, t+N-1),\Gamma)$. 
As we mentioned above, the  subset $S\subset(\al,\al+1)$ of the  zeros  of its determinant 
  is finite. Therefore, (\ref{seq}) implies  that each measure $\mu_j, 1\leq j<N,$ may only be concentrated at  $\{\al\}$ and on $S$, while the support of $\mu_N$ belongs to $\{\al,\al+1\}\cup S$. We may therefore write:
  $$
  d\mu_j=\sum_{s\in S}c_{s,j}\delta_{s}+a_j\delta_{\al}, \quad 1\leq j\leq N-1; $$$$ d\mu_N=\sum_{s\in S}c_{s,N}\delta_{s}+a_{N}\delta_{\al}+a_{N+1}\delta_{\al+1}.  $$
  This and (\ref{mu-decomp}) proves part (i) of the lemma, where $\nu$ and $\mu_j$ are defined in (\ref{rep}).
  
Finally, part (ii) easily follows from (\ref{seq}).
\end{proof}

\section{Uniqueness sets for lacunary polynomials}
In this section we will prove Theorem \ref{t2}. Clearly, if  $S(t_1, \ldots, t_N)$ is a uniqueness set for $P(M)$, then so is $-S(t_1, \ldots , t_N)$, since  $P(-x)\in P(M)$ whenever $P(x)\in P(M)$. 
Therefore, it suffices to prove that $S(t_1, \ldots , t_N)$ is a uniqueness set for every space $P(M),\#M=N.$

Assume a polynomial $P\in P(M)$ vanishes on $S(t_1, \ldots , t_N)$.
If $P$ is even or odd, we have $P(t_k) = 0, 1\le k \le N$, and by the Descartes' rule of signs we deduce that $P \equiv 0$. Thus, we can assume that $P\not\equiv0$ is neither even nor odd and derive a contradiction from there.

Consider the polynomials
$$P_e(x) = \sum_{m_j\in M, 2\mid m_j} c_j x^{m_j} = \frac{1}{2}(P(x) + P(-x))$$
and
$$P_o(x) = \sum_{m_j\in M, 2\nmid m_j} c_j x^{m_j} = \frac{1}{2}(P(x) - P(-x)).$$
If one of them is identically zero, then $P$ is even or odd and we are done. Let $M$ have $K$ even elements and $N-K$ odd elements. Then $P_e$ has at most $K-1$ positive roots and $P_o$ has at most $N-K-1$ positive roots by the Descartes' rule of signs. We are going to show that $P_e$ and $P_o$ together have at least $N-1$ positive roots thus getting the contradiction we need.

\begin{figure}[ht]
\center\includegraphics[height=8cm]{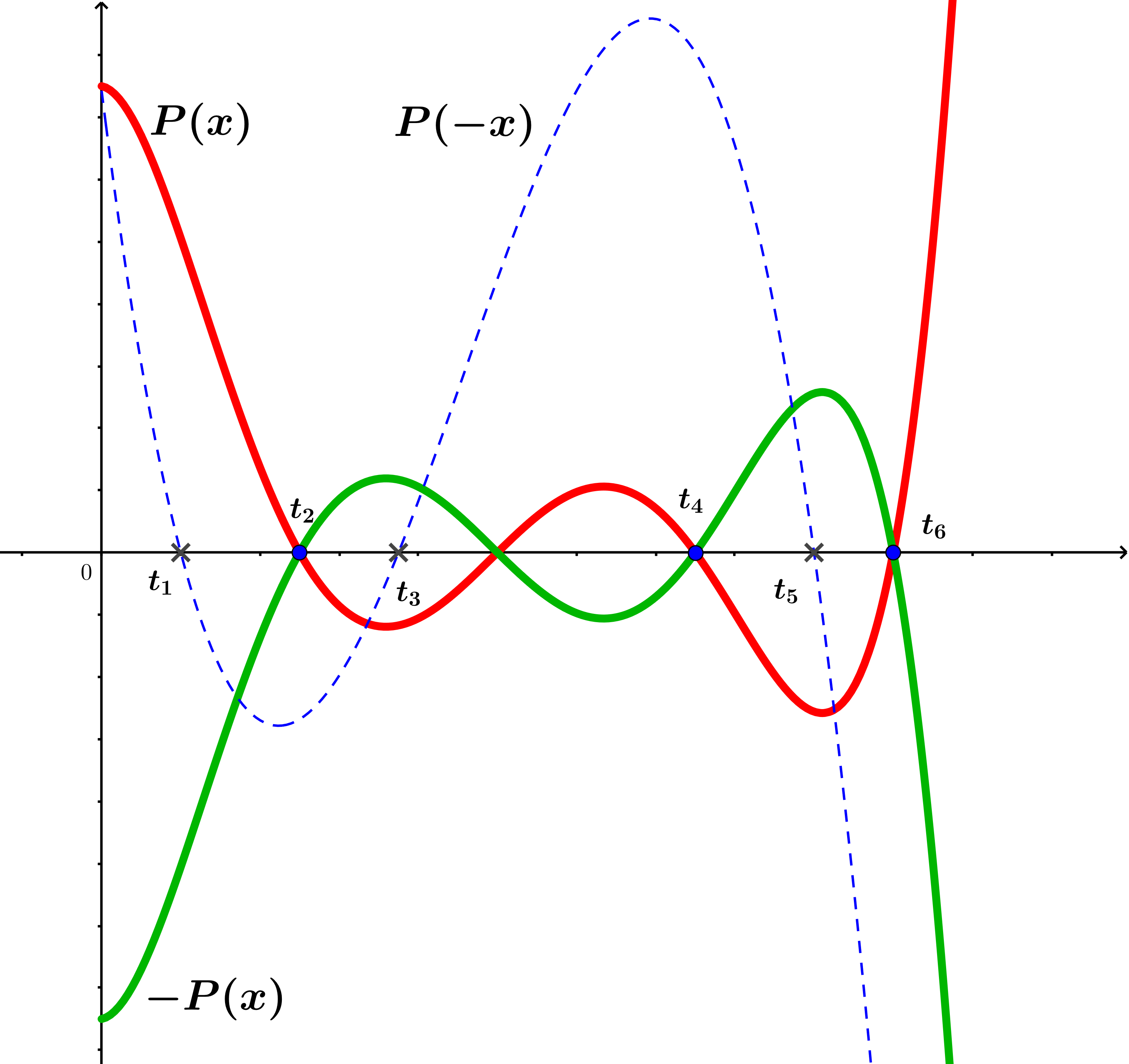}
\caption{$P(-x)$ has many intersections with either $P(x)$ or $-P(x)$}
\label{img1}
\end{figure}

Let us consider the graphs of  $P(x)$, $-P(x)$ and $P(-x)$, see Figure \ref{img1}. Since we assumed that $P$ is neither even nor odd, these are three different polynomials. For simplicity we first cover the case when $P(x)$ and $P(-x)$ do not have common positive zeroes. We indicate $t_k$ with odd indices by crosses.

By assumption each cross except the first one and the last one is separated from the other crosses by the zeroes of $P(x)$. That is, it is contained in a connected component bounded by the pieces of the curves $y = P(x)$ and $y = -P(x)$.  Thus, to get from the cross number $m$ to the cross number $m+1$ we have to exit the component containing the first and enter the next one, giving us at least two intersections of the curve $y = P(-x)$ with curves $y = P(x)$ and $y = -P(x)$. Additionally, if $N$ is even, then we also have to exit the last connected component as well, since there must be at least one more zero of $P(x)$ after the last cross. In total we will always have at least $N-1$ intersections, that is $P_e$ and $P_o$ together have at least $N-1$ positive roots as we wanted.

Now, we indicate the necessary changes in the case when $P(x)$ and $P(-x)$ have common positive roots. If we have two crosses which are not zeroes of $P(x)$ but between them there is a zero of $P(x)$, then the curve $y = P(-x)$ can go directly from the connected component of the first cross to the connected component of the second cross through this zero. But if $P(x_0) = P(-x_0) = 0$ then $x_0$ is a zero for both $P_e$ and $P_o$, thus we anyway get two zeroes.

It remains to consider the case when we have a cross which is also a zero of $P(x)$. Assume that crosses from the number $m$ to $m+l$ are zeroes of $P(x)$ and crosses number $m-1$ and $m+l+1$ are not (or there are no crosses with these indices). Then each of these $l+1$ zeroes are both zeroes for $P_e$ and $P_o$, thus giving us two intersections. Finally, since the $m+l$'th cross is separated from $m+l+1$'st by at least one more zero of $P(x)$ we have to enter the connected component corresponding to this zero and the same between $m$'th and $m-1$'st zero, thus giving us the same number of intersections as in the case when $P(x)$ and $P(-x)$ did not have common zeroes.

\section{Completeness of $E(\Z,\Gamma)$ in $L^2(a,b)$}

Part (i) of Theorem \ref{t1} follows from
\begin{theorem}\label{lemma_p} Given any finite set $\Gamma\subset\N_0$, the system
$E(\Z,\Gamma)$ is complete in $L^2(a,b)$ if and only if $b-a\leq \#\Gamma$.
\end{theorem}

\begin{proof}

(i) Assume $b-a\leq N:=\#\Gamma$. It is then a simple consequence of Lemma~\ref{l2} that $E(\Z,\Gamma)$
 is complete in $L^2(a,b)$. Indeed, if the system is not complete then there exists  non-trivial $f\in L^2(a, b)$ which is orthogonal to our system. Therefore, the measure $f\, dx$ is also orthogonal to the system, but it can not be a sum of delta measures unless $f$ is identically zero.
 
(ii) Assume that $b-a>N$.  We have to  prove that $E(\Z,\Gamma)$ is not complete in $L^2(a,b)$, i.e. that there is a non-trivial function $F\in L^2(a,b)$ such that
\begin{equation}\label{orth_f}
\int\limits_{a}^{b} t^\gamma e^{-2\pi i n t}\, F(t)\, dt=0,\quad \text{ for every    } \gamma\in\Gamma, n\in\Z.
\end{equation}
The existence of such a function follows essentially from elementary linear algebra.

We have  $b=a+N+\delta$, for some $\delta>0,$ and may assume that $\delta<1.$ 
Write $F$ in the form 
$$F(t)=\sum_{j=0}^NF_j(t-j),\quad t\in(a,a+N+\delta),$$where $F_j(t):=F(t+j){\bf 1}_{(a,a+1)}(t)$ vanish outside $(a,a+1)$ for $j=0,...,N-1$, and $f_N$ vanishes outside $(a,a+\delta)$. Here ${\bf 1}_{(a,a+1)}$ is the characteristic  function of $(a,a+1)$.
Clearly, to prove (\ref{orth_f}) it suffices to find $N+1$ non-trivial functions $F_j$ as above satisfying for a.e. $t\in (a,a+1)$ the system of $N$ equations  $$\sum_{j=0}^N (t+j)^\gamma F_j(t)=0,\quad \mbox{for all } \gamma\in\Gamma, \ t\in(a,a+1). $$ Rewrite this system in the matrix form 
$$V(t) \cdot  (F_0(t),...,F_{N-1}(t))^T=-((t+N)^{\gamma_1},...,(t+N)^{\gamma_N})^T\cdot F_N(t),\quad \Gamma=\{\gamma_1,...,\gamma_N\},$$
 where $V(t):=V((t, t+1, \dots, t+(N - 1);\Gamma)$ is a  generalized Vandermonde matrix defined above, whose determinant has only finite number of real zeroes. Therefore, there is an interval $I\subset(a,a+\delta)$ where $V(t)$ is invertible and satisfies $$\sup_{t\in I}\sup_{{\bf x}\in \R^N, \, \|x\|=1}\| V^{-1}(t)\cdot {\bf x}\|<\infty.$$
  
  Now, one can simply choose $F_N(t):={\bf 1}_{I}(t)$ and set  $$(F_0(t),...,F_{N-1}(t))^T:=-V^{-1}(t)\cdot((t+N)^{\gamma_1},...,(t+N)^{\gamma_N})^T\cdot{\bf 1}_{I}(t).$$
  \end{proof}

\begin{remark}
One can  check that the above result on completeness  of  $E(\Z, \Gamma)$ in $L^2(a,b)$ remain true for the space $L^p(a,b)$, $1 \le p < \infty.$
\end{remark}
\section{Completeness of $E(\Z,\Gamma)$ in $C([-a,a])$} 

\begin{theorem}\label{p}
$E(\Z,\Gamma)$ is complete in $C([-a,a])$ if and only if $a<r(\Gamma).$
\end{theorem}

Clearly, this theorem implies $CR_C(\Z,\Gamma)=r(\Gamma).$

\begin{proof}
 {\bf 1}. 
 Suppose $a\geq r(\Gamma)$. We have to check that the system is not complete in $C([-a,a])$. Clearly, it suffices  to produce a bounded measure $\mu$ on $[-r(\Gamma),r(\Gamma)]$ which satisfies (\ref{orth}).

Set $\mathbb{O}:=\#\Gamma_{odd}, \Ee:=\#\Gamma_{even}$ and 
\begin{equation}\label{f}
    f(x) = \begin{cases} \sin (\pi x) + 
 \sum\limits_{k=1}^{\mathbb{O}} \alpha_k \sin \left(\left(2k+1\right)\pi x\right) , \text{   if   }  \mathbb{O} < \Ee, \\[.2cm]
1+  \sum\limits_{k=1}^{\Ee} \alpha_k \cos \left(2\pi k x\right) , \text{   if   } \mathbb{O} \ge \Ee,\\[.2cm]
   \end{cases}
\end{equation}
where $\{\alpha_k\} \subset \R$.

\begin{lemma}\label{l4}
There exist numbers $\alpha_k$ in \eqref{f}
such that $f$  satisfies
\begin{equation} \label{interpol}
f^{(\gamma)}(n) = 0, \quad \gamma \in \Gamma,\,  n \in \N.\quad 
\end{equation}
\end{lemma}

It is easy to check that $f$ in (\ref{f}) is the Fourier-Stieltjes transform of a measure supported by $[-r(\Gamma),r(\Gamma)]$. One may therefore easily see that Lemma~\ref{l4} proves the necessity in part (i) of Theorem~\ref{p}.

\begin{proof}[Proof of Lemma \ref{l4}]
Consider the case $\Ee \le \mathbb{O}$.

We  wish to find $\alpha_k$ so that the function 
$$
f(x) = 1 + \alpha_1 \cos(2 \pi x) + \dots + \alpha_\Ee \cos(2 \pi \Ee x)
$$satisfies (\ref{interpol}).

It is clear that every {\it odd} derivative of  $f$ vanishes on $\Z$.
Therefore, it suffices  to find the coefficients so that   $f^{(\gamma)}$ vanishes  on $\Z$ for every $\gamma\in \Gamma_{even}$ (in particular, for $\gamma=0$). This is equivalent to saying that  the coefficients must satisfy the following  system of $\mathbb{E}$ linear equations:$$\gamma=0: \ \ \alpha_1+\dots +\alpha_\Ee=-1$$ and
$$\gamma\in \Gamma_{even}, \gamma\ne0: \ \ 
(2\pi)^{\gamma}\alpha_1 + (4\pi)^{\gamma}\alpha_2 \dots + (2\pi \Ee)^{\gamma}\alpha_\Ee = 0.
$$
 This system has a {\it unique non-trivial solution} by  Proposition~\ref{pos_matr_prop}.
  
The case $\Ee>\mathbb{O}$ is  similar, and we leave  the proof to the reader.
 \end{proof}

We return now to the proof of Theorem~\ref{p}.

{\bf 2}. 
Assume $a<r(\Gamma).$ We have to show that $E(\Z,\Gamma)$ is complete in $C([-a,a])$, i.e. that the only measure $\mu$ on $[-a,a]$ which satisfies (\ref{orth}) is trivial.

 We will consider the case $\Ee\le \mathbb{O}$, i.e. $r(\Gamma) = \Ee.$ Clearly, we may assume that  $\Ee= \mathbb{O}$ and so $\Ee=N/2$, where $N:=\#\Gamma$ is an even number. Also, to avoid trivial remarks, we assume that $N\geq 4$.

Assume that $\mu$ is concentrated on $[-a,a]$ and satisfies (\ref{orth}).
By (\ref{rep}) and Lemma~\ref{l2}, since $\mu(\{\pm N/2\})=0$, we have 
$$
d\mu=\sum_{s\in S}d\mu_s+d\nu= \sum_{s \in S} \sum_{j=1}^{N} c_{s,j} \delta_{s+j-1}+\sum_{j=2}^N a_j\delta_{-N/2+j-1}, 
$$where  $S$ is a finite subset of $(-N/2,-N/2+1)$ and  the coefficients $c_{s,j}$  satisfy for every $s\in S$ the system of equations  (\ref{coef}). By part (ii) of Lemma \ref{c1}, this system has only trivial solutions $c_{s,j}=0, j=1,...,N,s\in S\setminus (1/2)\Z$, and so  
$$
\mu= \nu_1+\nu,\quad  d\nu_1:=\sum_{j=1}^Nc_j\delta_{-N/2+j-1/2}, 
$$where $\nu$ and $\nu_1$ both are orthogonal to $E(\Z,\Gamma)$.

Let us check that $\nu=0.$
It is more convenient to write   $\nu$ in the form
$$
\nu=\sum_{k=-N/2+1}^{N/2-1}b_k\delta_k, \quad b_k:=a_{N/2+k+1}.
$$
Then clearly,  \eqref{orth} is equivalent to the system of $N-1$ equations: 
$$
\sum_{k=-N/2+1}^{N/2-1} k^\gamma b_k=0, \quad \text{ for every } \gamma \in \Gamma.
$$This is equivalent to the following systems:
$$
 \sum_{k=0}^{N/2-1} k^\gamma (b_{-k}+b_k)=0, \ \gamma\in\Gamma_{even},\ \ \sum_{k=1}^{N/2-1} k^\gamma (b_{-k}-b_k)=0, \ \gamma\in\Gamma_{odd}.
$$
One may now use Proposition \ref{pos_matr_prop} to deduce that  $b_{-k} + b_k= b_{-k} - b_k =0$, for every $k$, thus $b_k = b_{-k} = 0$ for every $k$, that is $\nu=0$.
Similarly, one may check that $\nu_1=0$, and so $\mu=0.$

The proof of the case $\mathbb{O} < \mathbb{E}$ is similar is left to the reader. 
\end{proof}

\begin{remark}
One can prove that for $a \in [r(\Gamma), \#\Gamma/2]$, the deficiency of $E(\Z, \Gamma)$ in $C([-a,a])$ is always finite.
\end{remark}

\section{Frame Property of $E(\Z,\Gamma)$}
The frame property of $E(\Z,\Gamma)$ in $L^2(a,b)$  is closely connected with the completeness property of $E(\Z,\Gamma)$  
in $C([a,b])$:

\begin{theorem}\label{ts}
Assume $a<b$ and $\epsilon>0.$

{\rm (i)} If $E(\Z,\Gamma)$ is complete in $C([a,b])$, then  $E(\Z,\Gamma)$ is a frame in $L^2(a,b).$ 

{\rm (ii)} If $E(\Z,\Gamma)$ is not complete in $C([a,b])$, then  $E(\Z,\Gamma)$ is not a frame in $L^2(a-\epsilon,b+\epsilon).$ 
\end{theorem}

Observe that to finish the proof of Theorem 1, it remains to show that $FR(\Z,\Gamma)=r(\Gamma).$ This follows easily from Theorem \ref{p} and  Theorem \ref{ts}.

\begin{proof}[Proof of Theorem \ref{ts}] (i) Assume that the system $E(\Z,\Gamma)$ is complete in $C([a,b])$. We have to show that it is a frame in $L^2(a,b)$.

Recall that  $E(\Z,\Gamma)$ is a frame in $L^2(a,b)$ if there are positive constants $A,B$ such that
\begin{equation}\label{eframe}
A\|F\|_2^2 \leq  \sum_{n\in\Z}\sum_{\gamma\in\Gamma}|\langle F,t^\gamma e^{2\pi i nt}\rangle|^2\leq B\|F\|_2^2,\quad \mbox{for every } F\in L^2(a,b).
\end{equation}Using the Fourier transform, this is equivalent to the condition\begin{equation}\label{eee}
A\|f\|_2^2 \leq  \sum_{n\in\Z}\sum_{\gamma\in\Gamma}|f^{(\gamma)}(n)|^2
\leq B\|f\|_2^2,\quad \end{equation}where $f$ is the inverse Fourier transform of $F$.

It is standard to check that the right hand-side inequality in (\ref{eframe}) (and in (\ref{eee})) holds for every interval $(a,b)$, see e.g. \cite{ou1}, Lecture 2. So, we only prove the left hand-side inequality.

By Theorem \ref{t1}, $E(\Z,\Gamma)$ is not complete, and so is not a frame in $L^2(a,b)$ when $b-a>N:=\#\Gamma.$ Therefore, in what follows we may assume that $a + k - 1 < b\leq a+k$, for some $k\in\N, k\leq N$.

Write
\begin{equation}\label{fff}
F(t)=\sum_{j=0}^{k-1}F_j(t-j),\quad F_j(t):=F(t+j)\cdot {\bf 1}_{(a,a+1)}(t).
\end{equation}Then we have
$$
 \langle F,t^\gamma e^{2\pi i nt}\rangle=\int\limits_a^{a+1}e^{2\pi i nt}\left(\sum_{j=0}^{k-1}(t+j)^\gamma F_j(t)\right)\,dt.
$$Hence,
$$
\sum_{n\in\Z}|\langle F,t^\gamma e^{2\pi i nt}\rangle|^2=
\|\sum_{j=0}^{k-1}(t+j)^\gamma F_j(t)\|_2^2.
$$
We see that the left hand-side inequality in (\ref{eframe}) is equivalent to 
\begin{equation}\label{f2}
\|V_k(t)\cdot (F_0(t),\dots,F_{k-1}(t))^T\|_2^2\geq A\|F\|_2^2,
\end{equation}where $$V_k(t):=V(t,\dots,t+ k-1;\Gamma)^T$$ denotes the $k\times N$ matrix which consists of the first $k$ columns of  $V(t,\dots, t+N-1;\Gamma)$, and we set $$\|(G_1,\dots,G_k)^T\|_2^2:=\|G_1\|_2^2+\dots+\|G_k\|_2^2.
$$

Let us first consider the case $b=a+k$.
Since $E(\Z,\Gamma)$ is complete in $C([a,b])$, there is no measure $\mu$ on $[a,b]$ orthogonal to this system. Then, since any measure of the form 
$$
d\mu=\sum_{j=0}^{k-1}x_j\delta_{t+j},\quad (x_1,\dots,x_k)\in\R^k\setminus\{\mathbf{0}\}, \ t\in[a,a+1],
$$is not orthogonal to $E(\Z,\Gamma)$, we see that
$V_k(t)\cdot {\bf x}^T\ne{\bf 0},$ for every ${\bf x}\in \R^k\setminus\{\mathbf{0}\}$ and $t\in[a,b]$. Therefore,
there is a constant $A$ such that $$\|V_k(t)\cdot {\bf x}^T\|^2\geq A\|{\bf x}\|^2,\quad t\in[a,a+1],$$ which implies  (\ref{f2}).

Now, let us assume that $b=a+k-1+\delta,$ where $0<\delta<1$. Then the function $F_{k-1}$ in (\ref{fff}) satisfies $F_{k-1}(t)=0, \delta<t<1.$ Similarly to above, for every  vectors ${\bf x}\in \R^k$ and ${\bf y}\in \R^{k-1}$ we have
$$
\|V_k(t)\cdot {\bf x}\|\geq A_1\|{\bf x}\|, \ t\in[a,a+\delta], \ \|V_{k-1}(t)\cdot {\bf y}\|\geq A_2\|{\bf y}\|, \ t\in[a+\delta,a+1],
$$from which (\ref{f2}) follows.

(ii) Assume that the system $E(\Z,\Gamma)$ is not complete in $C([a,b])$. We have to show that it is not a frame in $L^2(a-\epsilon,b+\epsilon)$, for every $\epsilon>0$. We may assume that $0<\epsilon<1/2.$ 

Let $g$ be the inverse Fourier transform of a measure $\mu$ on $[a,b]$ that is orthogonal to the system. Then $g^{(\gamma)}$ vanishes on $\Z$, for every $\gamma\in\Gamma$. 

Choose any $r$, $0<r<\epsilon,$ and consider the function 
$$
f(x):=g(x)\varphi(x),\quad \varphi(x):=\frac{\sin (\pi r x)}{\pi r x}.
$$
Then, clearly, $f$ is the (inverse) Fourier transform of an absolutely continuous measure on $(a-r,b+r)\subset(a-\epsilon,b+\epsilon)$, and
\begin{equation}\label{llll}\|f\|_2>C>0, \quad \mbox{where } C \mbox{  does not depend on } \epsilon.\end{equation} 

We will need
\begin{lemma}\label{lll}
There is a constant $C$ such that 
\begin{equation}\label{uu}\sum_{n\in\Z}|\varphi^{(j)}(n)|_2^2\leq C^j r^j,\quad j\in\N.\end{equation}
\end{lemma}

The proof of the lemma  follows from two observations: 

(i)  $\varphi$ is the  Fourier transform of ${\bf 1}_{(-r/2,r/2)}(t)/r$, and so $\varphi^{(j)}$ is the Fourier transform of $$(-2\pi i t)^j{\bf 1}_{(-r/2,r/2)}(t)/r.$$It easily follows that $\|\varphi^{(j)}\|_2^2\leq Cr^j,j\in\N.$

(ii) The sum in (\ref{uu}) is equal to the norm $\|\varphi^{(j)}\|_2^2$.

\medskip

Using (\ref{uu}), since $g^{(\gamma)},\gamma\in\Gamma,$ vanishes on $\Z$ and the functions $g^{(j)},j\in\N,$ are bounded on $\R$, one can easily check that $$
\sum_{n\in\Z}\sum_{\gamma\in\Gamma}|f^{(\gamma)}(n)|^2=\sum_{n\in\Z}\sum_{\gamma\in\Gamma}|(g\varphi)^{(\gamma)}(n)|^2\leq Cr,
$$for some $C$. This and (\ref{llll}) imply that the left hand-side inequality in (\ref{eee}) is not true for all small enough values of $r.$
\end{proof}

\begin{remark}
 Observe that by Theorem~{\rm \ref{lemma_p}}, $E(\Z,\Gamma)$ is not complete in $C([a,b])$ whenever $b-a>N:=\#\Gamma.$
Let us state  two results on  the completeness  of $E(\Z,\Gamma)$ in $C([a,b])$ when $a\geq0$:

{\rm (i)} Using part {\rm(i)} of Lemma~{\rm \ref{c1}} and Lemma~{\rm\ref{l2}}, one may check that $E(\Z,\Gamma)$ is complete in $C([a,b])$ whenever $b-a<N$ and if $a > 0$ then we don't need the assumption $0\in \Gamma$.

{\rm (ii)} One may also prove that $E(\Z,\Gamma)$ is complete in $C([a,a+N])$ if and only if $a\not\in\N_0.$
\end{remark}


\begin{remark}
 Let us come back to the exponential systems $E(\Z,\Gamma(n))$ defined in the beginning of  Section~{\rm{2}}. Here we present a simple example which illustrates   that such  systems  may have  strikingly different completeness properties in $L^2$-spaces and $C$-spaces.

Let $f(x)=\sin(\pi x/2)$. Then $f^{(2k)}(2n)=f^{(2k+1)}(2n+1)=0$, for every $k\in\N_0,n\in\Z$. Then, since $f$ is the inverse Fourier transform of $(\delta_{1/4}-\delta_{-1/4})/2i$,  the system
$$
\{t^{2k}e^{4\pi i nt}: k\in\N_0,n\in\Z\} \, \bigcup \, \{t^{2k+1} e^{2\pi i (2k+1)t}: k\in\N_0,n\in\Z\}
$$
is not complete in $C([-1/4,1/4])$. On the other hand, one may check that it is complete in $L^2(I)$ on every finite interval $I\subset\R$. 
\end{remark}

\section{Acknowledgements}
The authors want to thank Fedor Petrov and Pavel Zatitskiy for valuable discussions about this paper.

Aleksei Kulikov was supported by Grant 275113 of the Research Council of Norway, by BSF Grant 2020019, ISF Grant 1288/21,
and by The Raymond and Beverly Sackler Post-Doctoral Scholarship.

\noindent Aleksei Kulikov\\
Norwegian University of Science and Technology,\\ Department of Mathematical Sciences\\
NO-7491 Trondheim, Norway\\
\break
 Tel Aviv University,\\ School of Mathematical Sciences,\\ Tel Aviv 69978, Israel\\
lyosha.kulikov@mail.ru

\bigskip 

\noindent Alexander Ulanovskii\\
University of Stavanger, Department of Mathematics and Physics,\\
4036 Stavanger, Norway,\\
alexander.ulanovskii@uis.no

\bigskip

\noindent Ilya Zlotnikov\\
University of Stavanger, Department of Mathematics and Physics,\\
4036 Stavanger, Norway,\\
ilia.k.zlotnikov@uis.no

\end{document}